\documentclass[12pt]{amsart}

\usepackage{mathrsfs}

\usepackage{amsmath,amsthm,amssymb}
\font\teneufm=eufm10
\font\seveneufm=eufm7
\font\fiveeufm=eufm5
\newfam\frakturfam

\textfont\frakturfam=\teneufm
\scriptfont\frakturfam=\seveneufm
\scriptscriptfont\frakturfam=\fiveeufm

\newtheorem{pr}{Proposition}

\newtheorem{lm}{Lemma}
\newtheorem{theor}{Theorem}
\newtheorem{co}{Corollary}

\newtheorem{prob}{Problem}
\def\bee{\begin{eqnarray}}
\def\bes{\begin{eqnarray*}}
\def\eee{\end{eqnarray}}
\def\ees{\end{eqnarray*}}

\def\a{\alpha}
\def\b{\beta}

\def\Proof{{\sl Proof.}\ }

\title{Identities of the left-symmetric Witt algebras}

\begin{document}
\date{}
\maketitle

\begin{center}
{\bf Daniyar Kozybaev}\footnote{Supported by an MES grant 0755/GF of Kazakhstan; Eurasian National University,
 Astana, Kazakhstan,
e-mail: {\em kozybayev@gmail.com}}\\

and\\

{\bf Ualbai Umirbaev}\footnote{Supported by an MES grant 0755/GF of Kazakhstan; Eurasian National University,
 Astana, Kazakhstan and
 Wayne State University,
Detroit, MI 48202, USA,
e-mail: {\em umirbaev@math.wayne.edu}}
\end{center}

\begin{abstract} Let $P_n=k[x_1,x_2,\ldots,x_n]$ be the polynomial algebra over a field $k$ of characteristic zero in the variables $x_1,x_2,\ldots,x_n$ and $\mathscr{L}_n$ be the left-symmetric Witt algebra of all derivations of $P_n$ \cite{Burde06}. We describe all right operator identities of $\mathscr{L}_n$ and prove that the set of all algebras $\mathscr{L}_n$, where $n\geq 1$, generates the variety of all left-symmetric algebras. We also describe a class of general (not only right operator) identities for $\mathscr{L}_n$.
\end{abstract}

\noindent {\bf Mathematics Subject Classification (2010):} Primary 16S32,
17D25; Secondary 16R10, 17B01,17A50.
\noindent

{\bf Key words:} algebras of derivations, left-symmetric algebras, identities, free algebras, multiplication algebras.

\section{Introduction}

\hspace*{\parindent}

Let $k$ be an arbitrary field of characteristic zero  and $P_n=k[x_1,x_2,\ldots,x_n]$ be the polynomial algebra over $k$ in the variables $x_1,x_2,\ldots,x_n$. Let $W_n$ be the Witt algebra of index $n$, i.e., Lie algebra of all derivations of the polynomial algebra $P_n$.
The set of elements $u\partial_i$, where $u=x_1^{s_1}\ldots x_n^{s_n}\in P_n$ is an arbitrary monomial, $\partial_i=\frac{\partial}{\partial x_i}$, and $1\leq i\leq n$, forms a linear basis for $W_n$. Denote by  $\mathscr{L}_n$ an algebra with the ground space $W_n$ and the product $\cdot$  defined by
\bee\label{Df1}
a\partial_i\cdot b\partial_j= (a\partial_i(b))\partial_j,   \ \ \ a,b\in P_n.
\eee
 It is easy to check \cite{Burde06} that
$\mathscr{L}_n$ satisfies the left-symmetric identity
\bee\label{Df2}
(xy)z-x(yz)=(yx)z-y(xz).
\eee
This means that the associator $(x,y,z):=(xy)z-x(yz)$ is symmetric with respect to two left arguments, i.e.,
\bes
(x,y,z)=(y,x,z).
\ees
The variety of
left-symmetric algebras is Lie-admissible, i.e., each
left-symmetric algebra $\mathscr{L}$ with the operation $[x,y]:=xy-yx$ is a
Lie algebra. Obviously, the commutator algebra of $\mathscr{L}_n$ coincides with the Witt  algebra $W_n$ of all derivations of $P_n$. The algebra $\mathscr{L}_n$ is called the {\em left-symmetric Witt algebra of index $n$}.

The left-symmetric Witt algebra $\mathscr{L}_n$ of all derivations of the polynomial algebra $P_n$ plays an important role in the study of locally nilpotent derivations of $P_n$ and the Jacobian Conjecture \cite{UUTG}. The set of all left nilpotent elements of $\mathscr{L}_n$ coincides with the set of all locally nilpotent derivations of the polynomial algebra $P_n$ and right multiplication operator $R_D$, where $D\in \mathscr{L}_n$, is nilpotent if and only if the Jacobian matrix of $D$ (see, definition in \cite{UUTG}) is nilpotent. These facts already lead to necessity of purely algebraic study of $\mathscr{L}_n$. First of all, it is important to get algebraic descriptions of left and right nilpotent elements of $\mathscr{L}_n$.

Identities of $\mathscr{L}_n$ are especially important in order to treat the mentioned above  problems of affine algebraic geometry \cite{UUTG}.  Notice that
 $\mathscr{L}_1$ is a Novikov algebra, i.e., it also satisfies the identity
\bee\label{Df3}
(xy)z=(xz)y.
\eee
Moreover, the identities (\ref{Df2}) and (\ref{Df3}) form  a basis of identities for $\mathscr{L}_1$ in characteristic zero \cite{MLU11}. Recall that the question on a basis of identities for $W_1$ is a well-known open problem \cite{Razmyslov94}. A basis for  $\mathrm{Z}$-graded identities of the Lie algebra $W_1$ was recently given in \cite{FKK}.

The right multiplication algebra of $\mathscr{L}_n$, i.e., the associative algebra generated by all right multiplication operators $R_x, x\in \mathscr{L}_n$, is isomorphic to the matrix algebra $M_n(P_n)$ \cite[Lemma 4]{UUTG}. This result allows to make some parallel between the matrix algebras $M_n(k)$ in the case of associative algebras and the left-symmetric Witt algebras $\mathscr{L}_n$ in the case of left-symmetric algebras. In Section 2, we describe all right operator identities of the left-symmetric Witt algebra $\mathscr{L}_n$ as well as its subalgebras of triangular derivations $T(\mathscr{L}_n)$ and strongly triangular derivations $ST(\mathscr{L}_n)$. Recall that minimal right operator identities of $\mathscr{L}_n$ are considered in \cite{Dzhuma00}.

These results give that
\bes
\mathrm{Var}(\mathscr{L}_1)\varsubsetneqq \mathrm{Var}(\mathscr{L}_2)  \varsubsetneqq \ldots \varsubsetneqq \mathrm{Var}(\mathscr{L}_n)\varsubsetneqq \ldots,
\ees
where $\mathrm{Var}(\mathscr{L}_n)$ is the variety of algebras generated by $\mathscr{L}_n$.

In Section 4, we prove that the set of all left-symmetric Witt algebras $\mathscr{L}_n$, where $n\geq 1$, generates the variety of all left symmetric algebras. This means that every polynomial identity satisfied by all algebras $\mathscr{L}_n$ is a consequence of (\ref{Df2}). In fact, we prove more stronger result (Theorem \ref{Dt1}), which says that the algebra of derivations $ST(\mathscr{L}_n)$ does not satisfy any nontrivial left-symmetric identity of degree $\leq n$.

Unfortunately,  left operator identities of $\mathscr{L}_n$ are more important in applications and more difficult to describe. It is known that every Lie identity for $W_n$ is equivalent to a left operator identity for $\mathscr{L}_n$ \cite[Lemma 7]{UUTG}.
 Identities of the Lie algebra $W_n$ are studied in \cite{Razmyslov94}. In fact, Razmyslov's methods of proving identities for $W_n$ applied directly to $\mathscr{L}_n$ give more wider class of identities for $\mathscr{L}_n$. These identities are described in Section 5.

In Section 3, we give some auxiliary results on a basis of free left-symmetric algebras and some of them can be found in \cite{KMLU,Segal}.

\section{Algebras of triangular derivations}

\hspace*{\parindent}

Every element $D$ of the left-symmetric Witt algebra  $\mathscr{L}_n$ can be represented uniquely in the form
\bes
D_F=f_1\partial_1+f_2\partial_2+\ldots+f_n\partial_n,
\ees
where $F=(f_1,f_2,\ldots,f_n)'$ is an n-tuple (column) of elements of the polynomial algebra $P_n=k[x_1,x_2,\ldots,x_n]$. Let $J(F)=(\partial_j(f_i))_{1\leq i,j\leq n}$  be the Jacobian matrix of $F$ considered as a polynomial mapping of the space $k^n$. The Jacobian matrix $J(F)$ is also called the {\em Jacobian matrix of the derivation $D_F$}. Thereby, every derivation  $D\in \mathscr{L}_n$ has the Jacobian matrix $J(D)=J(D_F)=J(F)$.

Consider the grading
\bes
P_n=A_0\oplus A_1\oplus A_2\oplus\ldots\oplus A_s\oplus\ldots
\ees
of the polynomial algebra $P_n$, where $A_i$ is the space of homogeneous elements of degree $i\geq 0$. The left-symmetric algebra $\mathscr{L}_n$ has a natural grading
\bee\label{Xf1}
\mathscr{L}_n=L_{-1}\oplus L_{0}\oplus L_{1}\oplus \ldots \oplus L_{s}\oplus\ldots,
\eee
where $L_{i}$ is the space generated by all elements of the form $a\partial_j$ with $a\in A_{i+1}$ and $1\leq j\leq n$. Elements of $L_{s}$ are called {\em homogeneous derivations} of $P_n$ of degree $s$.

We have $L_{-1}=k\partial_1+\ldots+k\partial_n$ and $L_{0}$ is a subalgebra of $\mathscr{L}_n$ isomorphic to the matrix algebra $M_n(k)$. Notice that $L_{-1}\cdot L_n\subset L_{n-1}$ for all $n\geq 0$. If $0\neq f\in L_n$ and $n\geq 0$, then, obviously, there exists $\partial_s\in L_{-1}$ such that $\partial_s\cdot f\neq 0$. Following \cite{Razmyslov94}, this property will be referred as the {\em left transitivity} of the grading (\ref{Xf1}).

The element
\bes
D_X=x_1\partial_1+x_2\partial_2+\ldots+x_n\partial_n
\ees
 is the identity element of the matrix algebra $L_{0}$ and is the right identity element of $\mathscr{L}_n$. The left-symmetric algebra
$\mathscr{L}_n$ has no identity element.

Let $A$ be an arbitrary left-symmetric algebra. Denote by $\mathrm{Hom}_k(A,A)$ the associative algebra of all $k$-linear transformations of the vector space $A$.  For any $x\in A$ denote by $L_x : A\rightarrow A (a\mapsto xa)$ and $R_x : A\rightarrow A (a\mapsto ax)$ the operators of left and right multiplication by $x$, respectively. It follows from (\ref{Df2}) that
\bes
L_{[x,y]}=[L_x,L_y],\ \ \
R_{xy}=R_yR_x+[L_x,R_y].
\ees
Denote by $M(A)$ the subalgebra of $\mathrm{Hom}_k(A,A)$ (with an identity) generated by all $R_x,L_x$, where $x\in A$. The algebra $M(A)$ is called the {\em multiplication} algebra of $A$. The subalgebra $R(A)$ of $M(A)$ (with identity) generated by all $R_x$, where $x\in A$, is called the {\em right multiplication} algebra of $A$. Similarly, the subalgebra $L(A)$ of $M(A)$ (with identity) generated by all $L_x$, where $x\in A$, is called the {\em left multiplication} algebra of $A$.

For shortness of terminology, we say that $f=f(z_1,\ldots,z_m)$ is an {\em associative}  polynomial if  $f$ is an element of the free associative algebra $\mathrm{Ass}\langle z_1,\ldots,z_m\rangle$  freely generated by $z_1,\ldots,z_m$. If $f=f(z_1,\ldots,z_m)\in \mathrm{Ass}\langle z_1,\ldots,z_m\rangle$, then denote by $f^*$ the image of $f$ under the standard involution $*$ of $\mathrm{Ass}\langle z_1,\ldots,z_m\rangle$, i.e., $z_i^*=z_i$ for all $i$. We also say that $f=f(z_1,\ldots,z_m)$ is a {\em left-symmetric}  polynomial if  $f$ is an element of the free left-symmetric algebra $LS\langle z_1,\ldots,z_m\rangle$.

An identity $g=0$ of a left-symmetric algebra $A$ is called a {\em right operator identity} if $g=f(R_{x_1},\ldots,R_{x_m})x$ for some associative polynomial $f=f(z_1,\ldots,z_m)$. This means $f(R_{a_1},\ldots,R_{a_m})=0$ in $R(A)$ for all $a_1,\ldots,a_m\in A$. Notice that $f=0$ is not an identity for $R(A)$ (in general). Left operator identities can be defined similarly.

The right multiplication algebra $R(\mathscr{L}_n)$ of the left-symmetric Witt algebra $\mathscr{L}_n$ is isomorphic to the matrix algebra $M_n(P_n)$ \cite{UUTG}. This isomorphism
\bes
\theta : R(\mathscr{L}_n) \longrightarrow M_n(P_n)
\ees
is uniquely defined by $\theta(R_D)=J(D)$ for all $D\in \mathscr{L}_n$.

So the algebra
$R(\mathscr{L}_n)$ satisfies all identities of the matrix algebra $M_n(k)$. Obviously, every identity of $R(\mathscr{L}_n)$ determines an identity (in fact, a set of identities) for $\mathscr{L}_n$.
The isomorphism $\theta$ allows to carry out the study of $\mathscr{L}_n$ parallel to the study of $M_n(k)$. In this section we describe all right operator identities of $\mathscr{L}_n$ and its two important  subalgebras.

The next proposition gives a complete description of all  right operator identities of $\mathscr{L}_n$.

\begin{pr}\label{Dp1} Let $f=f(z_1,\ldots,z_m)$ be an associative polynomial. Then
$f(z_1,\ldots,z_m)=0$ is an identity for $M_n(k)$ if and only if $f(R_{y_1},\ldots,R_{y_m})y=0$ is a right operator identity for $\mathscr{L}_n$.
\end{pr}
\Proof Recall that $L_0\simeq M_n(k)$ is an associative algebra. We have
\bes
\partial_s \cdot (x_t\partial_r\cdot x_i\partial_j)=(\partial_s \cdot x_t\partial_r)\cdot x_i\partial_j
\ees
for all $s,t,r,i,j$. This means that $L_{-1}$ is a right $L_0$-module.

Transformation of matrices is an involution of $M_n(k)$. Consequently, $f=0$ is an identity for $M_n(k)$ if and only if $f^*=0$ is an identity for $M_n(k)$. Therefore, if $f=0$ is not an identity for $M_n(k)$, then there exist $a_1,\ldots,a_m\in L_0$ such that
$0\neq f^*(a_1,\ldots,a_m)\in L_0$. By the left transitivity of the grading (\ref{Xf1}), there exists $\partial\in L_{-1}$ such that $\partial(f^*(a_1,\ldots,a_m))\neq 0$. Since $L_{-1}$ is a right $L_0$-module, we get
\bes
\partial(f^*(a_1,\ldots,a_m))=f(R_{a_1},\ldots,R_{a_m})\partial\neq 0.
\ees
If $f=0$ is an identity for $M_n(k)$, then $f=0$ is an identity for $R(\mathscr{L}_n)$ since $R(\mathscr{L}_n)\simeq M_n(P_n)$. Consequently, $f(R_{a_1},\ldots,R_{a_m})=0$ for all $a_1,\ldots,a_m\in \mathscr{L}_n$. This means $f(R_{a_1},\ldots,R_{a_m})a=0$ for all $a\in \mathscr{L}_n$. $\Box$

According to the Amitsur--Levitzki theorem \cite{AL}, $M_n(k)$ satisfies the standard identity
\bes
s_{2n}(y_1,y_2,\ldots,y_{2n})=\sum_{\sigma\in S_{2n}} \mathrm{sgn}(\sigma) y_{\sigma(1)}y_{\sigma(2)}\ldots y_{\sigma(2n)}=0
\ees
of degree $2n$, where $S_{2n}$ is the symmetric group on degree $2n$  and $\mathrm{sgn}(\sigma)$ is the sign of the permutation $\sigma$. Moreover, $s_{2n}$ is an identity of $M_n(k)$ of the minimal degree \cite{AL}.

\begin{co}\label{Dc1} \cite{Dzhuma00} The left-symmetric Witt algebra
 $\mathscr{L}_n$ satisfies the identity
\bes
s_{2n}(R_{y_1},R_{y_2},\ldots,R_{y_{2n}})y_{2n+1}=0
\ees
of degree $2n+1$ and this is a right operator identity of $\mathscr{L}_n$ of the minimal degree.
\end{co}

 If $n=1$, this corollary gives the Novikov identity (\ref{Df3}).

Let $T(\mathscr{L}_n)$ be the subspace of $\mathscr{L}_n$ generated by all derivations of the form
\bes
x_i^{s_i}x_{i+1}^{s_{i+1}}\ldots x_n^{s_n}\partial_i,
\ees
where $s_i,s_{i+1},\ldots,s_n\in \mathbb{Z}_+$, $\mathbb{Z}_+$ is the set of all nonnegative integers, and $1\leq i\leq n$. It is easy to check that $T(\mathscr{L}_n)$ is a subalgebra of $\mathscr{L}_n$. We call $T(\mathscr{L}_n)$ the algebra of all {\em triangular} derivations of $P_n$. In fact, $T(\mathscr{L}_n)$ is the set of all derivations $D$ such that $J(D)$ is an upper triangular matrix, i.e., $\theta(R_D)=J(D)\in T_n(P_n)$. Recall \cite{Maltsev} that   the identity
\bes
[x_1,y_1][x_2,y_2]\ldots [x_n,y_n]=0
\ees
is a basis for all identities of the algebra $T_n(k)$ of all triangular matrices.

\begin{co}\label{Dc2} Let $f=f(z_1,\ldots,z_m)$ be an arbitrary element of the free associative algebra $\mathrm{Ass}\langle z_1,\ldots,z_m\rangle$ and let $C$ be the commutator ideal of $\mathrm{Ass}\langle z_1,\ldots,z_m\rangle$. Then the following conditions are equivalent.

(i) $f(z_1,\ldots,z_m)=0$ is an identity for $T_n(k)$;

(ii) $f\in C^n$;

(iii) $f(R_{y_1},\ldots,R_{y_m})y=0$ is an identity for $T(\mathscr{L}_n)$.
\end{co}
\Proof Notice that (i) and (ii) are equivalent \cite{Maltsev}. This gives that $f=0$ is an identity for $T_n(k)$ if and only if $f^*=0$ is an identity for $T_n(k)$. Notice that $L_{-1}\subset T(\mathscr{L}_n)$ and $T_n(k)=L_0\cap T(\mathscr{L}_n)\subset T(\mathscr{L}_n)$. Using these facts, we can repeat the proof of Proposition \ref{Dp1} for this case without any changes.  $\Box$

Denote by $ST(\mathscr{L}_n)$ the subspace of $\mathscr{L}_n$ generated by all derivations of the form
\bes
x_{i+1}^{s_{i+1}}\ldots x_n^{s_n}\partial_i,
\ees
where $s_{i+1},\ldots,s_n\in \mathbb{Z}_+$ and $1\leq i\leq n$. Then $ST(\mathscr{L}_n)$ is a subalgebra of $T(\mathscr{L}_n)$. We call $ST(\mathscr{L}_n)$ the algebra of all {\em strongly triangular} derivations of $P_n$. Notice that $ST(\mathscr{L}_n)$ is the set of all derivations $D$ such that $J(D)$ is an upper strongly triangular matrix, i.e., $J(D)\in ST_n(P_n)$. The algebra $ST_n(k)$ of all strongly triangular matrices satisfies the identity
\bes
y_1y_2\ldots y_n=0.
\ees
It is easy to check that $ST_n(k)$ does not satisfy any nontrivial associative identity of degree less than $n$.

Repeating the proofs of Proposition \ref{Dp1} and Corollary \ref{Dc2}, we get the next corollary.
\begin{co}\label{Dc3} Every right operator identity of the algebra $ST(\mathscr{L}_n)$ of all strongly triangular derivations is a corollary of the identity
\bee\label{Df4}
R_{y_1}R_{y_2}\ldots R_{y_n}z=0.
\eee
\end{co}

So $ST(\mathscr{L}_n)$ is right nilpotent. It is not difficult to show that $ST(\mathscr{L}_n)$ is locally nilpotent but not nilpotent. In fact, the identity of right nilpotency does not give too much in the case of left-symmetric algebras. An example of a two-dimensional right nilpotent but not nilpotent Novikov algebra is given in \cite{Zelmanov87}. More strange examples of Novikov algebras in small dimensions can be found in \cite{BG}. The picture might be totally different for subalgebras of $\mathscr{L}_n$.

\begin{prob}\label{Dpr1} Let $\mathscr{L}$ be a subalgebra of $\mathscr{L}_n$ satisfying the identity $R_y^nz=0$. Is $\mathscr{L}$ locally nilpotent?
\end{prob}

Identities of the form $L_y^mz=0$ have stronger corollaries \cite{Filippov01} in the case of left-symmetric algebras and related to the study of Engel Lie algebras  \cite{Zelmanov88}.

\section{Reduced words}

\hspace*{\parindent}

Let $Y=\{y_1,y_2,\ldots,y_n\}$ be an alphabet of $n$ symbols. Denote by
$Y^*$ the monoid of all nonassociative words on $Y$ \cite{KBKA}. If $u\in Y^*$, then denote by
$d(u)$ the length of $u$. Every nonassociative word $u$ of length $\geq 2$ can be uniquely written as $u=u_1u_2$, where $d(u_1), d(u_2)<d(u)$.
If $u\in Y^*$ and $d(u)\geq 2$, then we consider only nontrivial decompositions of $u$, i.e.,
if $u=u_1u_2$, then $d(u_1),d(u_2)\geq 1$.

For any $u\in Y^*$ denote by $[u]$ the associative word obtained from $u$ by skipping all parentheses.

Put $y_1<y_2<\ldots <y_n$. Let $u$ and $v$ be arbitrary elements
of $Y^*$. We say that $u<v$ if $d(u)<d(v)$. If $d(u)=d(v)\geq 2$,
$u=u_1u_2$, and $v=v_1v_2$, then $u<v$ if either $u_1<v_1$ or
$u_1=v_1$ and $u_2<v_2$.

Put $\bar{n}=\{1,2,\ldots,n\}$. Every mapping $\sigma: \bar{n}\rightarrow \bar{n}$ will be identified with the endomorphism $\sigma: Y^*\rightarrow Y^*$ of the monoid $Y^*$ defined by $\sigma(y_i)=y_{\sigma(i)}$.

\begin{lm}\label{Dl1}
Let $u,v\in Y^*$ be arbitrary two words such that $u>v$, $[u]=y_{i_1}y_{i_2}\ldots y_{i_m}$,  $[v]=y_{j_1}y_{j_2}\ldots y_{j_r}$, and let $\sigma: \bar{n}\rightarrow \bar{n}$ be a mapping such that
$\sigma(i_1)>\ldots >\sigma(i_m)>\sigma(j_1)>\ldots >\sigma(j_r)$. Then $\sigma(u)>\sigma(v)$.
\end{lm}
\Proof Assume that the statement of the proposition is not true and let $u,v\in Y^*$ be a counterexample with the minimal $d(u)+d(v)$. If $d(u)>d(v)$, then $\sigma(u)>\sigma(v)$ since
$d(\sigma(u))>d(\sigma(v))$. Consequently, we may assume that $d(u)=d(v)$. If $d(u)=d(v)=1$ then $u=y_{i_1}$, $v=y_{j_1}$. Therefore,  $\sigma(u)=y_{\sigma(i_1)}$, $\sigma(v)=y_{\sigma(j_1)}$, and $\sigma(u)>\sigma(v)$ since $\sigma(i_1)>\sigma(j_1)$.

If $d(u)=d(v)\geq 2$, then $u=u_1u_2$ and $v=v_1v_2$. If $u_1>v_1$, then $u_1,v_1$ satisfies all the conditions of the proposition and $d(u_1)+d(v_1)<d(u)+d(v)$. By the choice of $u,v$ we get
$\sigma(u_1)>\sigma(v_1)$. Consequently, $\sigma(u)>\sigma(v)$. If $u_1=v_1$, then $u_2>v_2$. Similarly, we get $\sigma(u_2)>\sigma(v_2)$ and $\sigma(u)>\sigma(v)$.
$\Box$

A word $w\in Y^*$ is called  {\em reduced} \cite{Segal}  if it does not contain any subword of the form
$r(st)\in Y^*$, where
 $d(r),d(s),d(t)\geq 1$ and $r<s$. Denote by $V$ the set of all reduced words in the alphabet $Y$.

\begin{lm}\label{Dl2}
Let $w\in V$, $[w]=y_{i_1}y_{i_2}\ldots y_{i_m}$, and let $\sigma: \bar{n}\rightarrow \bar{n}$ be a mapping such that
$\sigma(i_1)>\ldots> \sigma(i_m)$. Then $\sigma(w)\in V$.
\end{lm}
\Proof   If $d(w)=1$, then there is nothing to prove. Suppose that $d(w)\geq 2$ and the statement of the proposition is true for all reduced words of length less than $d(w)$.
Let $w=uv$. Then $\sigma(u), \sigma(v)\in V$  by the induction proposition.  If $d(v)=1$, then , obviously, $\sigma(w)=\sigma(u)\sigma(v)\in V$. If $v=v_1v_2$, then
$\sigma(w)=\sigma(u)(\sigma(v_1)\sigma(v_2))$. We have $u\geq v_1$ since $w\in V$.
By Lemma \ref{Dl1}, we get $\sigma(u)\geq \sigma(v_1)$. This implies that $\sigma(w)\in V$ since $\sigma(u),\sigma(v)\in V$. $\Box$

Let $A=LS\langle Y\rangle$ be the free
left-symmetric algebra with free set of generators $Y$.  Every nonassociative word in the alphabet $Y$
represents a certain element of $A$ and for any $u\in Y^*$ we denote by $u$ the
element of $A$ defined by $u$.

According to \cite{Segal},  the set of all reduced words forms a linear basis for $A$:
every nonzero element $g$ of $A$ can be
uniquely represented as
\bee\label{Df5}
g=\a_1 w_1+\a_2 w_2+\ldots +\a_m w_m,
\eee
where $w_i\in V$, $0\neq \a_i\in k$ for all $i$, and $w_1<w_2<\ldots <w_m$.

Denote by $\overline{g}$ the {\em lowest word} $w_1$ of $g$.
The element $\a_1w_1$ is also called the {\em lowest term} of $g$.

\begin{lm}\label{Dl3}
Let $w\in Y^*$ be an arbitrary nonassociative word.
Then $\overline{w}\geq w$ in the algebra $A=LS\langle Y\rangle$ and the equality holds if and only if $w\in V$.
\end{lm}
\Proof
Assume that the statement of the lemma is true for all words of length less than $t$, where $t\geq 2$. Assume also that $w$ is the maximal word of length $t$ for which the statement of the lemma is not true. If $w\in V$, then
$\overline{w}=w$.  Hence $w=uv\notin V$. If
$u$  or $v$ is not reduced, then, by the assumptions above, $w$ is a linear
combination of words of the form $u'v'$ such that $u'v'>w$. By the choice of $w$, the statement of the lemma is true for all $u'v'$. Consequently, $w$ is a linear
combination of reduced words $w'$ such that $w'>w$.

Suppose that  both $u$ and $v$ are reduced. This means $v=v_1v_2$ and
$u<v_1$ since $w$ is not reduced. By (\ref{Df2}), we have
\bes
w=u(v_1v_2)=v_1(uv_2)+(uv_1)v_2-(v_1u)v_2
\ees
in the algebra $A$.
Notice that $v_1(uv_2),(uv_1)v_2,(v_1u)v_2>w$. By the choice of $w$, we again get that $w$ is a linear
combination of reduced words $w'$ such that $w'>w$. $\Box$

\begin{lm}\label{Dl4}
Every reduced word $w\in V$ can be uniquely represented in the form
\bee\label{Df6}
w=L_{w_1}L_{w_2}\ldots L_{w_m}x_i,
\eee
where $w_j\in V$ for all $j$ and $w_1\geq w_2\geq \ldots \geq w_m$.
\end{lm}
\Proof If $w=uv$, then $u$ and
$v$ are reduced words. Leading an induction on $d(w)$, we may assume that
$v=L_{v_1}L_{v_2}\ldots L_{v_s}x_i$, where $v_j\in V$ and $v_1\geq
v_2\geq \ldots \geq v_s$. We have $v=v_1v'$ where
$v'=L_{v_2}\ldots L_{v_s}x_i$. Note that $u\geq v_1$
since $w=uv=v(v_1v')$ is reduced. Consequently,
$w=L_uL_{v_1}L_{v_2}\ldots L_{v_s}x_i$ and $u\geq v_1\geq v_2\geq \ldots
\geq v_s$. The
uniqueness  of this representation is obvious. $\Box$

\section{Generalized triangular derivations}

\hspace*{\parindent}

Let $\lambda=(\lambda_{12},\lambda_{13},\ldots,\lambda_{1n},\lambda_{23},\ldots,\lambda_{n-1 n})$ be a tuple of independent commutative variables $\lambda_{ij}$, where $1\leq i<j\leq n$.
Denote by $k[\lambda]$ the polynomial algebra over $k$ in the variables $\lambda_{ij}$. Every element of $k[\lambda]$ can be written uniquely in the form
\bes
\lambda^s=\lambda_{12}^{s_{12}}\lambda_{13}^{s_{13}}\ldots \lambda_{1n}^{s_{1n}}\lambda_{23}^{s_{23}}\ldots \lambda_{n-1 n}^{s_{n-1 n}},
\ees
where $s=(s_{12},s_{13},\ldots,s_{1n},s_{23},\ldots,s_{n-1 n})\in \mathbb{Z}_+^{r}$ and $r=\frac{(n-1)n}{2}$. Denote by $\leq$ the lexicographic order on $\mathbb{Z}_+^{r}$. Put  also $\lambda^s<\lambda^t$ if $s<t$ in $\mathbb{Z}_+^{r}$. If $0\neq f\in k[\lambda]$ then denote by $\widetilde{f}$ the leading monomial of $f$.

Denote by $S$ the set of formal symbols
\bes
x_1^{f_1}\ldots x_n^{f_n},
\ees
where $f_1,f_2,\ldots,f_n\in k[\lambda]$, and consider $S$ as a semigroup with respect to the product
\bes
x_1^{f_1}\ldots x_n^{f_n}x_1^{g_1}\ldots x_n^{g_n}=x_1^{f_1+g_1}\ldots x_n^{f_n+g_n}.
\ees
Notice that $S$ just is a multiplicative version of the additive semigroup $k[\lambda]^n$. Put also $\deg_{x_i}(x_1^{f_1}\ldots x_n^{f_n})=f_i$ for all $i$.

Consider the set of formal symbols $u\partial_i$,
where $u\in S$ and $1\leq i\leq n$.
Let $\mathscr{L}$ be the free $k[\lambda]$-module generated by all $u\partial_i$. We turn $\mathscr{L}$ into $k[\lambda]$-algebra by
\bes
u\partial_i \circ v\partial_j=\deg_{x_i}(v) (uvx_i^{-1})
\partial_j.
\ees
\begin{lm}\label{Dl5}
$\mathscr{L}$ is a left-symmetric $k[\lambda]$-algebra.
\end{lm}
\Proof If $u,v,w\in S$, then
\bes
(u\partial_i,v\partial_j,w\partial_k)=(u\partial_i \circ v\partial_j)\circ w\partial_k-u\partial_i\circ(v\partial_j\circ w\partial_k)\\
=\deg_{x_i}(v)(uvx_i^{-1})\partial_j\circ w\partial_k-\deg_{x_j}(w)u\partial_i\circ (vwx_j^{-1})\partial_k\\
=\deg_{x_i}(v)\deg_{x_j}(w)(uvwx_i^{-1}x_j^{-1})\partial_k-\deg_{x_j}(w)\deg_{x_i}(vwx_j^{-1})(uvw x_i^{-1}x_j^{-1})\partial_k\\
=(\deg_{x_i}(v)\deg_{x_j}(w)-\deg_{x_j}(w)\deg_{x_i}(vwx_j^{-1}))(uvw x_i^{-1}x_j^{-1})\partial_k\\
=-\deg_{x_j}(w)(\deg_{x_i}(w)+\deg_{x_i}(x_j^{-1}))(uvw x_i^{-1}x_j^{-1})\partial_k.
\ees
Similarly, we get
\bes
(v\partial_j,u\partial_i,w\partial_k)=-\deg_{x_i}(w)(\deg_{x_j}(w)+\deg_{x_j}(x_i^{-1}))(uvw x_i^{-1}x_j^{-1})\partial_k.
\ees
If $i=j$, then these equalities give
\bes
(u\partial_i,v\partial_j,w\partial_k)=(v\partial_j,u\partial_i,w\partial_k),
\ees
i.e.,
the identity (\ref{Df2}) holds for $u\partial_i,v\partial_j,w\partial_k$. If $i\neq j$, then $\deg_{x_i}(x_j^{-1})=\deg_{x_j}(x_i^{-1})=0$ and (\ref{Df2}) holds again.
$\Box$

Denote by $\widehat{\mathscr{L}}_n$ the left-symmetric algebra of all derivations of $k[x_1^{\pm 1},x_2^{\pm 1},\ldots,x_n^{\pm 1}]$. The set of elements $u\partial_i$, where $u=x_1^{s_1}\ldots x_n^{s_n}$, $s_1,\ldots,s_n\in \mathbb{Z}$, $\partial_i=\frac{\partial}{\partial x_i}$, and $1\leq i\leq n$, forms a linear basis for $\widehat{\mathscr{L}}_n$. The product in $\widehat{\mathscr{L}}_n$ is defined by (\ref{Df1}).  The identity (\ref{Df2}) can be checked as in the proof of Lemma \ref{Dl5}. The subalgebras of triangular derivations $T(\widehat{\mathscr{L}}_n)$ and strong triangular derivations $ST(\widehat{\mathscr{L}}_n)$ of $\widehat{\mathscr{L}}_n$ can be defined similarly.

For any $s=(s_{12},s_{13},\ldots,s_{1n},s_{23},\ldots,s_{n-1 n})\in \mathbb{Z}^{r}$ we define a $k$-linear mapping
\bee\label{Df7}
\widehat{s}: \mathscr{L}\longrightarrow \widehat{\mathscr{L}}_n
\eee
by
\bes
\widehat{s}(f(\lambda)x_1^{f_1(\lambda)}\ldots x_n^{f_n(\lambda)}\partial_i)=
f(s)x_1^{f_1(s)}\ldots x_n^{f_n(s)}\partial_i
\ees
for all $f,f_1,\ldots,f_n\in k[\lambda]$ and $1\leq i\leq n$. Obviously, $\widehat{s}$ is defined correctly since $\mathscr{L}$ is a free $k[\lambda]$-module and $k[\lambda]$ is a free $k$-module.

\begin{lm}\label{Dl6}
$\widehat{s}$ is a homomorphism of $k$-algebras.
\end{lm}
\Proof We have
\bes
f(\lambda)x_1^{f_1(\lambda)}\ldots x_n^{f_n(\lambda)}\partial_i \circ g(\lambda)x_1^{g_1(\lambda)}\ldots x_n^{g_n(\lambda)}\partial_j\\
=f(\lambda)g(\lambda)g_i(\lambda)x_1^{f_1(\lambda)+g_1(\lambda)}\ldots x_i^{f_i(\lambda)+g_i(\lambda)-1}\ldots x_n^{f_n(\lambda)+g_n(\lambda)}\partial_j
\ees
in $\mathscr{L}$.
Obviously, we have
\bes
f(s)x_1^{f_1(s)}\ldots x_n^{f_n(s)}\partial_i \cdot g(s)x_1^{g_1(s)}\ldots x_n^{g_n(s)}\partial_j\\
=f(s)g(s)g_i(s)x_1^{f_1(s)+g_1(s)}\ldots x_i^{f_i(s)+g_i(s)-1}\ldots x_n^{f_n(s)+g_n(s)}\partial_j
\ees
in $\widehat{\mathscr{L}}_n$.
This means that $\widehat{s}$ is a homomorphism of $k$-algebras.
$\Box$

 Put
\bes
u_1=x_2^{\lambda_{12}}\ldots x_n^{\lambda_{1n}}, u_2=x_3^{\lambda_{23}}\ldots x_n^{\lambda_{2n}},\ldots, u_n=1.
\ees

Consider the elements
\bes
z_1=u_1\partial_1, z_2=u_2\partial_2,\ldots, z_n=u_n\partial_n,
\ees
of $\mathscr{L}$ and denote by $\mathscr{L}(\lambda)$ the $k$-subalgebra of $\mathscr{L}$ generated by $z_1,z_2,\ldots,z_n$.

For any $s=(s_{12},s_{13},\ldots,s_{1n},s_{23},\ldots,s_{n-1 n})\in \mathbb{Z}_+^{r}$,  consider the homomorphism $\widehat{s}$ defined by (\ref{Df7}). Denote the restriction of $\widehat{s}$ into $\mathscr{L}(\lambda)$ by the same symbol $\widehat{s}$.
Notice that $\widehat{s}(\mathscr{L}(\lambda))\subset ST(\mathscr{L}_n)$ since $\widehat{s}(z_i)\in ST(\mathscr{L}_n)$. Finally, we have a homomorphism
\bee\label{Df8}
\widehat{s} : \mathscr{L}(\lambda)\longrightarrow ST(\mathscr{L}_n).
\eee

As in Section 3, let $A=LS\langle Y\rangle$ be the free left-symmetric algebra generated by $Y=\{y_1,y_2,\ldots,y_n\}$.
Denote by
\bes
\chi : A \longrightarrow \mathscr{L}(\lambda)
\ees
the homomorphism of $k$-algebras such that $\chi(y_i)=z_i$ for all $i$.

\begin{lm}\label{Dl7} Let $w\in Y^*$. Then there exist unique polynomials $f^w, f_1^w, \ldots, f_n^w\in k[\lambda]$ and a number $r(w)(1\leq r(w)\leq n)$ such that
\bes
\chi(w)=f^w x_1^{f_1^w}\ldots x_n^{f_n^w}\partial_{r(w)}.
\ees
\end{lm}
\Proof If $w=y_i$ then $\chi(w)=z_i$. Hence $f^w=1, f_1^w=\ldots=f_i^w=0, f_{i+1}^w=\lambda_{i i+1},\ldots, f_n^w=\lambda_{i n}$, and $r(w)=i$. If $w=uv$, then
\bes
\chi(w)=\chi(u)\circ\chi(w)=f^u x_1^{f_1^u}\ldots x_n^{f_n^u}\partial_{r(u)}\circ f^v x_1^{f_1^v}\ldots x_n^{f_n^v}\partial_{r(v)}\\
=f^uf^v f_{r(u)}^v x_1^{f_1^u+f_1^v}\ldots x_{r(u)}^{f_{r(u)}^u+f_{r(u)}^v-1}\ldots x_n^{f_n^u+f_n^v}\partial_{r(v)}
\ees
by the definition of the product $\circ$ since $\deg_{x_{r(u)}}(x_1^{f_1^v}\ldots x_n^{f_n^v})=f_{r(u)}^v$. $\Box$

We use this lemma also as the definition of $f^w$, $f_i^w$, and $r(w)$. The proof of this lemma also implies the next corollaries.
\begin{co}\label{Dc4} Let $w\in Y^*$  and $w=uv$. Then
$f^w=f^uf^v f_{r(u)}^v$, $f_i^w={f_i^u+f_i^v}$ if $i\neq r(u)$, $f_{r(u)}^w={f_{r(u)}^u+f_{r(u)}^v}-1$, and $r(w)=r(v)$.
\end{co}
\begin{co}\label{Dc5} Let $w\in Y^*$ and $[w]=y_{i_1}\ldots y_{i_m}$. Then $r(w)=i_m$ and
\bes
\chi(w)=f^w u_{i_1}u_{i_2}\ldots u_{i_m} x_{i_1}^{-1}x_{i_2}^{-1}\ldots x_{i_{m-1}}^{-1}
\partial_{i_m}.
\ees
\end{co}
\Proof Formula $r(w)=r(v)$ from Corollary \ref{Dc4} immediately gives $r(w)=i_m$.

If $d(w)=1$ and $w=y_{i_1}$, then $\chi(w)=z_{i_1}=u_{i_1}\partial_{i_1}$. Let $w=uv$,  $[u]=y_{i_1}\ldots y_{i_m}$, and $[v]=y_{j_1}\ldots y_{j_r}$. Leading an induction on $d(w)$, we get
\bes
\chi(w)=\chi(u)\circ\chi(v)= f^u u_{i_1}\ldots u_{i_m} x_{i_1}^{-1}\ldots x_{i_{m-1}}^{-1}
\partial_{i_m} \circ f^v u_{j_1}\ldots u_{j_r} x_{j_1}^{-1}\ldots x_{j_{r-1}}^{-1}
\partial_{j_r} \\
=f^w u_{i_1}\ldots u_{i_m}u_{j_1}\ldots u_{j_r} x_{i_1}^{-1}\ldots x_{i_{m-1}}^{-1}x_{j_1}^{-1}\ldots x_{j_{r-1}}^{-1}x_{i_m}^{-1}
\partial_{j_r}
\ees
since $f^w=f^uf^vf_{i_m}^v=
f^uf^v \deg_{x_{i_m}}(u_{j_1}\ldots u_{j_r} x_{j_1}^{-1}\ldots x_{j_{r-1}}^{-1})$.  This implies the statement of the corollary. $\Box$

By this corollary, $r(w)$ is the index of the rightmost symbol of $w$. We use this fact as a denotation.

Recall that a word $w\in Y^*$ is called {\em multilinear} if $d_{y_i}(w)\leq 1$ for all $i$, where $d_{y_i}(w)$ is the length of $w$ in $y_i$.

A word $w\in Y^*$ with $[w]=y_{i_1}y_{i_2}\ldots y_{i_m}$ is called an {\em $s$-word} if
$i_1,i_2,\ldots,i_{m-1}>i_m$. We say $w\in Y^*$ is {\em special} if   every subword of $w$ is an  $s$-word.

Notice that if $w$ is a reduced word, $[w]=y_{i_1}y_{i_2}\ldots y_{i_m}$, and $i_1>i_2>\ldots>i_m$, then $w$ is a special reduced word. Obviously, every special reduced word $w$ has a unique representation in the form (\ref{Df6}), where $w_1,w_2,\ldots,w_m$ are special reduced words.

\begin{lm}\label{Dl8} Let $v\in Y^*$ be a multilinear word. If $v$ is not special, then $\chi(v)=0$
\end{lm}
\Proof Let $u$ be a nonspecial subword of $v$ of the minimal length. Then $u$ is not an $s$-word. Let $[u]=y_{i_1}\ldots y_{i_m}$ and let $i_j$ be the smallest number between $i_1,i_2,\ldots,i_m$. Notice that all numbers
$i_1,i_2,\ldots,i_m$ are different since $v$ is multilinear and $j<m$ since $u$ is not an $s$-word. Consequently, $i_{j+1},i_{j+2},\ldots,i_m>i_j$. This means that
$\chi(y_{i_{j+1}})=z_{i_{j+1}},\chi(y_{i_{j+2}})=z_{i_{j+2}},\ldots,\chi(y_{i_m})=z_{i_m}$ do not depend on $x_{i_j}$.
Consequently, $\chi(y_{i_j})=z_{i_j}$ annihilates the product $\chi(y_{i_1})\ldots \chi(y_{i_m})$ for any placement of parentheses. $\Box$

Denote by $W$ the set of all multilinear special reduced words in the alphabet  $Y$.
If $w\in W$ and $[u]=y_{i_1}\ldots y_{i_m}$, then put $s(w)=\{i_1,i_2,\ldots,i_m\}$.

\begin{lm}\label{Dl9} Let $w\in W$ and $w=uv$. Then
\bes
f_{r(u)}^v=\sum_{i\in s(v)}  \lambda_{i r(u)}.
\ees
\end{lm}
\Proof Recall that $\lambda_{ij}=0$ if $i\geq j$. Consequently, only the indexes $i\in s(v)$ with $i<r(v)$ contribute to the sum in the formulation of the lemma.

Notice that $v$ does not contain $y_{r(u)}$ since $w$ is multilinear and $y_{r(u)}$ is the rightmost symbol of $u$. Consequently, we do not have  $z_{r(u)}$ in $\chi(v)$. This means that the degrees in $x_{r(u)}$ of all elements participating in the product $\chi(v)$ will be summarized. Notice that
the degree $x_{r(u)}$ in $z_i$ is $\lambda_{i r(u)}$ if $i<r(u)$ and $0$ otherwise. This gives the statement of the lemma.
$\Box$

\begin{lm}\label{Dl10} Let $w\in W$ represented in the form (\ref{Df6}).
 Then
\bes
\widetilde{f^w}= \lambda_{i r(w_1)}\lambda_{i r(w_2)}\ldots \lambda_{i r(w_m)}
\widetilde{f^{w_1}}\widetilde{f^{w_2}}\ldots \widetilde{f^{w_m}}.
\ees
\end{lm}
\Proof We have $w=w_1u$, where $u=L_{w_2}\ldots L_{w_m}y_i$. By Corollary \ref{Dc4},
$f^w=f^{w_1}f^u f_{r(w_1)}^u$. By Lemma \ref{Dl9}, we get
\bes
f_{r(w_1)}^u=\sum_{j\in s(u)}  \lambda_{j r(w_1)}.
\ees
Notice that $i$ is the least number of $s(u)$ and $s(w)$ since $u$ and $w$ are special. Consequently, $\widetilde{f_{r(w_1)}^u}=\lambda_{i r(w_1)}$ and
\bes
\widetilde{f^w}=\lambda_{i r(w_1)}\widetilde{f^{w_1}}\widetilde{f^u}.
\ees
Leading an induction on $d(w)$, we may assume that
\bes
\widetilde{f^u}=\lambda_{i r(w_2)}\ldots \lambda_{i r(w_m)}
\widetilde{f^{w_2}}\ldots \widetilde{f^{w_m}}.
\ees
This completes the proof of the lemma. $\Box$.

\begin{co}\label{Dc6} In the conditions of Lemma \ref{Dl10} the following statements are true:

(i) $\widetilde{f^w}$ is a product of $d(w)-1$ elements $\lambda_{st}$;

(ii) The least first index of all divisors $\lambda_{st}$ of $\widetilde{f^w}$ is $i$;

(iii) The set of all second indexes of all divisors $\lambda_{st}$ of $\widetilde{f^w}$ is equal to $s(w)\setminus \{i\}$;

(iv) Let $\lambda_{i j_1},\lambda_{i j_2},\ldots,  \lambda_{i j_p}$ be the set of all divisors $\lambda_{st}$ of $\widetilde{f^w}$ with the first index $i$. Then $p=m$ and
$\{j_1,j_2,\ldots,j_p\}=\{r(w_1),r(w_2),\ldots,r(w_m)\}$.
\end{co}
\Proof Every statement of the corollary is obvious or can be deduced from Lemma \ref{Dl10} by an easy induction on $d(w)$. $\Box$

\begin{lm}\label{Dl11} Every $w\in W$ is uniquely defined by $\widetilde{f^w}$.
\end{lm}
\Proof  We determine $w\in W$ in the form (\ref{Df6}). By Corollary \ref{Dc6}(ii), $r(w)=i$ is determined uniquely. By Corollary \ref{Dc6}(iv), the number $m$ and the set of elements $\{r(w_1),r(w_2),\ldots,r(w_m)\}$ are also determined  uniquely.

We introduce new orders $\succ$ and $\vdash$ on $s(w)$. If $p,q\in s(w)$, then put $q\succ p$ if  $\lambda_{pq}$ divides $\widetilde{f^w}$. Put also $q\vdash p$ if $q=q_0\succ q_1\succ\ldots \succ q_t=p$ for some $t\geq 1$. By Corollary \ref{Dc6}(iii), for any $q\in s(w)\setminus \{i\}$ there exists a unique element $p\in s(w)$ such that $q\succ p$. Lemma \ref{Dl10} implies that $q\vdash r(w_j)$ if and only if $q\in s(w_j)\setminus \{r(w_j)\}$. This fact uniquely determines $s(w_j)$ for all $j$. Notice that $\widetilde{f^{w_j}}$ is the product of all divisors $\lambda_{st}$ of $\widetilde{f^w}$ such that $t\vdash r(w_j)$. Leading an induction on $d(w)$, we may assume that $w_j$ is uniquely determined by $\widetilde{f^{w_j}}$.
Then $w$ is also determined uniquely since $w_1\geq w_2\geq\ldots\geq w_m$. $\Box$

\begin{theor}\label{Dt1} The left-symmetric algebra  of all strongly triangular derivations   $ST(\mathscr{L}_n)$ does not satisfy any nontrivial left-symmetric identity of degree less than or equal to $n$.
\end{theor}
\Proof Recall \cite{KBKA} that  every polynomial identity over a field of characteristic zero is equivalent to homogeneous multilinear polynomial identities. Since $ST(\mathscr{L}_{n-1})\subset ST(\mathscr{L}_n)$, it is sufficient to prove that every multilinear polynomial identity of $ST(\mathscr{L}_n)$ of degree $n$ is a corollary of (\ref{Df2}). Suppose that it is not true. Then there exists a nonzero multilinear element  $g=g(y_1,y_2,\ldots,y_n)$ of the free left-symmetric algebra $A=LS\langle Y\rangle$ which is a polynomial  identity for $ST(\mathscr{L}_n)$, i.e.,
\bes
g(a_1,a_2,\ldots,a_n)=0
\ees
 for all $a_1,a_2,\ldots,a_n\in \mathscr{L}_n$. We fix this element $g$ and assume that $g$ is written in the form (\ref{Df5}). Then $\overline{g}=w_1$

Let  $[\overline{g}]=[w_1]=y_{i_1}y_{i_2}\ldots y_{i_n}$. We have $\{i_1,i_2,\ldots,i_n\}=\{1,2,\ldots,n\}$ since $g$ is multilinear of degree $n$. Consider the mapping $\sigma : \bar{n} \rightarrow \bar{n}$ defined $\sigma(i_j)=n-j+1$ for all $j$. Notice that $\sigma(i_1)>\sigma(i_2)>\ldots >\sigma(i_n)$. Denote by $\sigma$ the automorphism of the monoid $Y^*$ and the automorphism of the free algebra $A$ defined by $\sigma$. We have
\bes
\sigma(g)=\a_1 \sigma(w_1)+\a_2 \sigma(w_2)+\ldots +\a_m \sigma(w_m).
\ees
By Lemma \ref{Dl1}, we get $\sigma(w_1)<\sigma(w_2)<\ldots<\sigma(w_m)$.  By Lemma \ref{Dl2}, $\sigma(w_1)$ is a reduced word. Then Lemma \ref{Dl3} gives that $\overline{\sigma(g)}=\sigma(w_1)$.

Obviously, $\sigma(w_1)$ is special since
$[\sigma(w_1)]=y_ny_{n-1}\ldots y_1$.

Consequently, by changing $g$ to $\sigma(g)$, we may assume that
 $\overline{g}=w_1\in W$. Let $g=h+h_1$, where $h$ is a linear combination of special reduced words and $h_1$ is a linear combination of reduced nonspecial words. We have $h\neq 0$ since
$\overline{g}\in W$. We also have $\chi(h_1)=0$ by Lemma \ref{Dl8}. Assume that
\bes
h=\b_1 v_1+\b_2 v_2+\ldots +\b_s v_s,
\ees
where $v_i\in W$, $0\neq \b_i\in k$, $s\geq 1$, and $1\leq i\leq s$. Notice that $r(v_i)=1$ for all $i$ since $v_i$ is special. By Corollary \ref{Dc5}, we get $\chi(v_i)=f^{v_i} u_1u_2\ldots u_n(x_2x_3\ldots x_n)^{-1}\partial_1$ for all $i$. Hence
\bes
\chi(g)=\chi(h)=(\b_1 f_{v_1}+\b_2 f_{v_2}+\ldots +\b_s f_{v_s})u_1u_2\ldots u_n(x_2x_3\ldots x_n)^{-1}\partial_1.
\ees
Put $f_g=\b_1 f_{v_1}+\b_2 f_{v_2}+\ldots +\b_s f_{v_s}$. By Lemma \ref{Dl9}, the elements $\widetilde{f_{v_i}}$ are different for different values of $i$. Consequently, $f_g\in k[\lambda]$ is nonzero. It is not difficult to show that there exists $s\in \mathbb{Z}_+^r$ such that $f_g(s)\neq 0$. Then the image of $g$ under the homomorphism
\bes
\widehat{s}\circ \chi : A\longrightarrow ST(\mathscr{L}_n)
\ees
is nonzero. Consequently, $g$ is not an identity for $ST(\mathscr{L}_n)$. $\Box$

\begin{co}\label{Dc7} The identity (\ref{Df4}) is a left-symmetric polynomial identity of $ST(\mathscr{L}_n)$ of the minimal degree.
\end{co}
\begin{co}\label{Dc8} The variety of algebras generated by all algebras of strongly triangular derivations $ST(\mathscr{L}_n)$, where $n\geq 1$, is the variety of all left-symmetric algebras.
\end{co}
\begin{co}\label{Dc9} The variety of algebras generated by all algebras of triangular derivations $T(\mathscr{L}_n)$, where $n\geq 1$, is the variety of all left-symmetric algebras.
\end{co}
\begin{co}\label{Dc10} The variety of algebras generated by all left-symmetric Witt algebras $\mathscr{L}_n$, where $n\geq 1$, is the variety of all left-symmetric algebras.
\end{co}

\section{General identities}

\hspace*{\parindent}

A very interesting class of identities for $W_n$ was discovered by Yu.P. Razmyslov \cite[Chapter 6]{Razmyslov94}. Of course, every identity for $W_n$ is an identity for $\mathscr{L}_n$.  Moreover, every identity for $W_n$ also gives a left operator identity for $\mathscr{L}_n$ \cite[Lemma 7]{UUTG}. But Razmyslov's method, applied directly to $\mathscr{L}_n$, gives more wider class of identities for $\mathscr{L}_n$. We describe this class of identities.

Consider the grading (\ref{Xf1}) of $\mathscr{L}_n$. Recall that $0\neq a\in L_i$ is called homogeneous of degree $i$. Put $|a|=i$ in this case. Choose a homogeneous basis
\bes
e_1,e_2,\ldots, e_s,\ldots
\ees
for $\mathscr{L}_n$ such that $|e_i|\leq |e_j|$ if $i<j$. Then $e_1,\ldots,e_n$ is a basis for $L_{-1}$ and $e_{n+1},\ldots,e_{n^2+n}$ is a basis for $L_0$, and so on.

For any positive integer $N$ put
\bes
e(N)=\sum_{i=1}^N |e_i|
\ees

Recall that a (nonassociative) polynomial $f=f(y_1,\ldots,y_N,z_1,\ldots,z_t)$ is called {\em skew-symmetric} with respect to $y_1,\ldots,y_N$ if
\bes
\sigma(f)=f(y_{\sigma(1)},\ldots,y_{\sigma(N)},z_1,\ldots,z_t)
=\mathrm{sgn}(\sigma)f(y_1,\ldots,y_N,z_1,\ldots,z_t)
\ees
for all $\sigma\in S_N$, where $S_N$ is the group of degree $N$.

\begin{pr}\label{Dp2} Let $N$ be an arbitrary positive integer and
 $f=f(y_1,\ldots,y_N,z_1,\ldots,z_t)$ be an arbitrary multilinear left-symmetric polynomial skew-symmetric with respect to $y_1,\ldots,y_N$. If $e(N)\geq t$, then $f=0$ is an identity for $\mathscr{L}_n$.
\end{pr}
\Proof The identity (\ref{Df2}) gives
\bes
\partial\cdot(u\cdot v)=(\partial\cdot u)\cdot v+u\cdot (\partial\cdot v)
\ees
for all $\partial\in L_{-1}$ and $u,v\in \mathscr{L}_n$ since $\mathscr{L}_n\cdot L_{-1}=0$.
This means that $L_{\partial}$ is a derivation of $\mathscr{L}_n$.

Let $g=g(z_1,\ldots,z_l)$ be an arbitrary multilinear left-symmetric polynomial.
If $g=0$ is not an identity for $\mathscr{L}_n$, then there exist homogeneous elements $v_1,\ldots,v_l\in \mathscr{L}_n$ such that $0\neq g(v_1,\ldots,v_l)\in L_{-1}$.
In fact, let
$v_1,\ldots,v_l\in \mathscr{L}_n$ be homogeneous elements such that $0\neq g(v_1,\ldots,v_l)\in L_i$ and $i\geq 0$. By the left transitivity of the grading (\ref{Xf1}), there exists $\partial\in L_{-1}$ such that $L_{\partial} g(v_1,\ldots,v_l)\neq 0$. We have
\bes
L_{\partial} g(v_1,\ldots,v_l)=\sum_{i=1}^l g(v_1,\ldots,L_{\partial} v_i,\ldots,v_l).
\ees
since $L_{\partial}$ is a derivation. Then there exists $i$ such that $0\neq  g(v_1,\ldots,L_{\partial} u_i,\ldots,v_l)\in L_{i-1}$.

Suppose that $f=0$ is not an identity for $\mathscr{L}_n$. Then there exist basis elements $e_{i_1},\ldots,e_{i_N},e_{j_1},\ldots,e_{j_t}\in \mathscr{L}_n$ such that $0\neq f(e_{i_1},\ldots,e_{i_N},e_{j_1},\ldots,e_{j_t})\in L_{-1}$.
If two of the elements $e_{i_1},\ldots,e_{i_N}$ are equal to each other, then $f(e_{i_1},\ldots,e_{i_N},e_{j_1},\ldots,e_{j_t})=0$ since $f$ is skew-symmetric with respect to $y_1,\ldots,y_N$. If all $e_{i_1},\ldots,e_{i_N}$ are different, then
\bes
|e_{i_1}|+\ldots+|e_{i_N}|\geq |e_1|+\ldots+|e_N|=e(N).
\ees
This implies that $f(e_{i_1},\ldots,e_{i_N},e_{j_1},\ldots,e_{j_t})\in L_{s}$ and $s\geq e(N)-t\geq 0$. $\Box$

Let's describe the identities of the minimal degree given by Proposition \ref{Dp2}.
Notice that $|e_1|+\ldots+|e_n|=-n$ and $|e_{n+1}|+\ldots+|e_{n^2+n}|=0$. The minimal $N$ for which $e(N)$ is nonnegative is $n^2+2n$ and $e(n^2+2n)=0$. Put $N=n^2+n$. Let $w$ be an arbitrary reduced word such that $[w]=y_1\ldots y_N$. Then
\bes
S_N^w=\sum_{\sigma\in S_N} \mathrm{sgn}(\sigma)w(y_{\sigma(1)},\ldots,y_{\sigma(N)})=0
\ees
is an identity for $\mathscr{L}_n$.

\end{document}